\documentclass[]{article}

\usepackage{amsmath,amsfonts,amssymb} 

\usepackage{graphicx}       
\usepackage{subfigure}
\usepackage{array}
\usepackage{vmargin}

\providecommand{\keywords}[1]{\textbf{\textit{Key words:}} #1}

\providecommand{\subclassname}[1]{\textbf{\textit{Subject Classifications:}} #1}

\providecommand{\acknowledgements}[1]{\textbf{\textit{Acknowledgements:}} #1}

\newcommand{\C}{\mathbb{C}} 
\newcommand{\R}{\mathbb{R}} 
\newcommand{\N}{\mathbb{N}} 

\DeclareSymbolFont{bbold}{U}{bbold}{m}{n}
\DeclareSymbolFontAlphabet{\mathbbold}{bbold}
\newcommand{\ind}{\mathbbold{1}}

\newcommand{\PP}{\mathbb{P}}
\newcommand{\EE}{\mathbb{E}}
\newtheorem{algo}{Algorithm}
\newtheorem{theorem}{Theorem}
\newtheorem{remark}{Remark}
\newtheorem{proposition}{Proposition}
\newtheorem{lemma}{Lemma}
\newcommand{\Proof}{\noindent {\em{Proof. }}}
\newcommand{\EProof}{\begin{flushright}$\Box$\end{flushright}}
\newtheorem{assumption}[theorem]{Assumption}

\usepackage{microtype} 

\usepackage[colorlinks=true,linkcolor=black,citecolor=black,urlcolor=black]{hyperref}
\usepackage{bookmark}
\makeatletter 
  \providecommand*{\toclevel@author}{999}
  \providecommand*{\toclevel@title}{0}
\makeatother

\begin{document}

\title{Central Limit Theorem for Adaptive Multilevel Splitting Estimators in an Idealized Setting}
\author{
Charles-Edouard Br\'ehier
\footnote{
Universit\'e Paris-Est, CERMICS (ENPC),  6-8-10 Avenue Blaise Pascal, Cit\'e Descartes, F-77455 Marne-la-Vall\'ee, France. 
INRIA Paris-Rocquencourt, Domaine de Voluceau - Rocquencourt, B.P. 105 - 78153 Le Chesnay, France.
e-mail~:~brehierc@cermics.enpc.fr, e-mail~:~charles-edouard.brehier@inria.fr
}
\and Ludovic Gouden\`ege
\footnote{
F\'ed\'eration de Math\'ematiques de l'\'Ecole Centrale Paris, CNRS, Grande voie des vignes, 92295 Ch\^atenay-Malabry, France.
e-mail~:~goudenege@math.cnrs.fr
}
\and Lo\"ic Tudela
\footnote{
Ensae ParisTech, 3 Avenue Pierre Larousse, 92240 Malakoff, France.
e-mail~:~loic.tudela@ensae-paristech.fr
}
}

\maketitle

\abstract{The Adaptive Multilevel Splitting algorithm is a very powerful and versatile iterative method to estimate the probability of rare events, based on an interacting particle systems.
In \cite{brehier-lelievre-rousset-2014}, in a so-called idealized setting, the authors prove that some associated estimators are unbiased, for each value of the size $n$ of the systems of replicas and of a resampling number $k$.\newline
Here we go beyond and prove these estimator's asymptotic normality when $n$ goes to $+\infty$, for any fixed value of $k$. The main ingredient is the asymptotic analysis of a functional equation on an appropriate characteristic function.
 Some numerical simulations illustrate the convergence to rely on Gaussian confidence intervals.
}


\subclassname{ 65C05; 65C35; 60F05}

\section{Introduction}

Many models from physics, chemistry or biology involve stochastic systems for different purposes: taking into account some uncertainty with respect to some data parameters, or to allow for dynamical phase transitions between different configurations of the system. This phenomenon often referred to as metastability is observed for instance when one uses a $d$-dimensional overdamped Langevin dynamics
$$dX_t=-\nabla V(X_t)dt+\sqrt{2\beta^{-1}}dW_t$$
associated with a potential function $V$ with several local minima. Here $W$ denotes a $d$-dimensional standard Wiener process. When the inverse temperature $\beta$ increases, the transitions become rare events (their probability decreases exponentially fast).



In this paper, we adopt a numerical point of view, and analyze a method which outperforms a pure Monte-Carlo method for a given computational effort in the small probability regime (in terms of relative error). Two important families of methods have been introduced in the $1950$s and extensively developed later in order to efficiently address this rare event estimation problem: importance sampling, and importance/multilevel splitting - see \cite{kahn-harris} and \cite{glass_et_al} for a more recent treatment. For more general references, we refer for instance to \cite{rubino-tuffin}.

The method we study in this work belongs is a multilevel splitting algorithms. The main advantage of this kind of methods is that they are non-intrusive: one does not need to modify the model to go beyond a pure Monte-Carlo approximation scheme. The present method has an additional feature: the so-called levels are computed adaptively. To explain more precisely the algorithm and its properties, from now on we only focus on a simpler setting, though paradigmatic of the rare event estimation problem.

Let $X$ be a real random variable, and $a$ some fixed, given threshold. We want to compute an approximation of the tail probability $p:=\mathbb{P}(X>a)$. The splitting strategy in the regime when $a$ becomes large consists in introducing the following decomposition of $p$ as a product of conditional probabilities:
$$\PP(X>a)=\PP(X>a_n | X>a_{n-1})\ldots \PP(X>a_{2} | X>a_1)\PP(X>a_1),$$
for some sequence of levels $a_1<\ldots<a_{n-1}<a_n=a$. The common interpretation of this formula is that the event that $X>a$ is split in $n$ (conditional) crossing probabilities for $X$, each much higher than $p$  and thus easier to approximate (each independently).

To optimize the variance, the levels must be chosen such that all the conditional probabilities are equal to $p^{1/n}$, with $n$ as large as possible. However, the \textit{a priori} knowledge of levels satisfying this condition is not available in practical cases.

Notice that in principle to apply this splitting strategy, one needs to know how to sample according to the conditional distributions appearing in the splitting formula. If this condition holds, we say that we are in an idealized setting.

Adaptive techniques, where the levels are computed on-the-fly have been introduced in the $2000$s in various contexts, under different names: Adaptive Multilevel Splitting (AMS) \cite{cerou-guyader-07a}, \cite{cerou-guyader-07b}, \cite{cerou-guyader-lelievre-pommier-11}, Subset simulation \cite{au_beck} and Nested sampling \cite{skilling}.

In this paper, we focus on the versions of AMS algorithms studied in \cite{brehier-lelievre-rousset-2014}. Two parameters are required: a number of (interacting) replicas $n$, and a fixed integer $k\in\left\{1,\ldots,n-1\right\}$, such that the proportion of replicas that are killed and resampled at each iteration is $k/n$. The version with $k=1$ has been studied in \cite{guyader-hengartner-matzner-lober-11}, and is also (in the idealized setting) a special case of the Adaptive Last Particle Algorithm of \cite{simonnet}.

A family of estimators $(\hat{p}^{n,k})_{n\geq 2,1\leq k\leq n-1}$ is introduced in \cite{brehier-lelievre-rousset-2014}. The main property established there is unbiasedness: $\mathbb{E}[\hat{p}^{n,k}]=p$. Moreover, the computational cost is analyzed in the regime when $k$ is fixed and $n\rightarrow +\infty$. Nevertheless comparisons between different values of $k$ is made assuming a non natural procedure: $N$ independent realizations of the algorithm are necessary to define an empirical estimator with realizations of $\hat{p}^{n,k}$. We remove this procedure by showing directly an asymptotic normality result for the estimator $\hat{p}^{n,k}$, allowing for the direct use of Gaussian asymptotic confidence intervals.

Other Central Limit Theorems for Adaptive Multilevel Splitting estimators (in different parameter regimes) have been obtained in \cite{cerou-del-moral-furon-guyader-12}, \cite{cerou-guyader-07a} and \cite{cerou-guyader-delmoral-malrieu-14}.

The main result of this paper is Theorem \ref{th}: if $k$ and $a$ are fixed, under the assumption that the cumulative distribution function of $X$ is continuous, when $n\rightarrow +\infty$ we have the convergence in law of $\sqrt{n}\bigl(\hat{p}^{n,k}-p\bigr)$ to a centered Gaussian random variable, with variance $-p^2\log(p)$ (independent on $k$).

The main novelty of the paper is the case $k>1$: indeed when $k=1$ the law of the estimator is explicitly known (it involves a Poisson random variable with parameter $-n\log(p)$) and the Central Limit Theorem can be derived by hand. More precisely, we prove the asymptotic normality of $\log(\hat{p}^{n,k})$, and conclude thanks to the delta-method. However when $k>1$, the key idea is to prove a functional equation, as introduced in \cite{brehier-lelievre-rousset-2014}, for the characteristic function of this random variable; the basic ingredient is a decomposition according to the first step of the algorithm. Thus one of the main messages of this paper is that the functional equation technique allows to prove several key properties of the AMS algorithm in the idealized setting: unbiasedness and asymptotic normality.

The paper is organized as follows. In Section \ref{sect:AMS}, we introduce the main objects: the idealized setting (Section \ref{subsect:setting}) and the AMS algorithm (Section \ref{subsect:algo}). Our main result (Theorem \ref{th}) is stated in Section \ref{subsect:CLT}. Section \ref{sect:proof} is devoted to the detailed proof of this result. Finally Section \ref{sect:num} contains a numerical illustration of the Theorem.

\section{Adaptive Multilevel Splitting Algorithms}\label{sect:AMS}

\subsection{Setting}\label{subsect:setting}

Let $X$ be some real random variable. We assume that $X>0$ almost surely.
We want to estimate the probability $p=\PP(X>a)$, where $a>0$ is some threshold. When $a$ goes to $+\infty$, $p$ goes to $0$ and we have to estimate a rare event.
More generally, we introduce the conditional probability for $0\leq x\leq a$
\begin{equation}
P(x)=\PP(X>a |X>x).
\end{equation}

We notice that $p=P(0)$ and that $P(a)=1$.

The following assumption is crucial:
\begin{assumption}\label{ass:cdf:c0}
Let $F$ denote the cumulative distribution function of $X$. We assume that $F$ is continuous.
\end{assumption}

\subsection{Algorithm}\label{subsect:algo}

The algorithm depends on two parameters, fixed once and for all:
\begin{itemize}
\item the number of replicas $n$;
\item the number $k\in\left\{1,\ldots,n-1\right\}$ of replicas that are resampled at each iteration.
\end{itemize}

The other necessary parameters are the initial condition $x$ and the stopping threshold $a$.

In the sequel, when we consider a random variable $X_{i}^{j}$, the subscript $i$ denotes the index in $\left\{1,\ldots,n\right\}$ of a particle, and the superscript $j$ denotes the iteration of the algorithm.

In the algorithm below and in the following, we use classical notations for $k$-th order statistics. For $Y=(Y_1,\ldots,Y_n)$ independent and identically distributed (i.i.d.) real valued random variables with continuous cumulative distribution function, there exists  almost surely a unique (random) permutation $\sigma$ of $\left\{1,\ldots,n\right\}$ such that $Y_{\sigma(1)}<\ldots<Y_{\sigma(n)}$. For any $k \in \{1, \ldots,n\}$, we then use the classical notation $Y_{(k)}=Y_{\sigma(k)}$ to denote the $k$-th order statistics of the sample $Y$.

\begin{algo}[Adaptive Multilevel Splitting]\label{algo:AMS}
~

\noindent
{\bf Initialization:}
Define $Z^{0}=x$.
Sample $n$ i.i.d. realizations $X_{1}^{0},\ldots,X_{n}^{0}$, with the law $\mathcal{L}(X | X>x)$.

Define $Z^{1}=X_{(k)}^{0}$, the $k$-th order statistics of the sample $X^{0}=(X_{1}^{0},\ldots,X_{n}^{0})$, and $\sigma^1$ the (a.s.) unique associated permutation: $X_{\sigma^1(1)}^{0}<\ldots<X_{\sigma^1(n)}^{0}$. 

Set $j=1$.

\noindent
{\bf Iterations (on $j\geq 1$):} While $Z^{j} <  a$:

\begin{itemize}
\item Conditionally on $Z^{j}$, sample $k$ new independent random variables $(Y_1^j,\ldots,Y_k^j)$, according to the law $\mathcal{L}(X | X>Z^{j})$.

\item Set
$$
X_{i}^{j}=\begin{cases}Y_{(\sigma^j)^{-1}(i)}^{j} \quad \text{if } (\sigma^j)^{-1}(i)\leq k\\ X_{i}^{j-1} \quad \text{if } (\sigma^j)^{-1}(i)>k. \end{cases}
$$

In other words, the particle with index $i$ is killed and resampled according to the law $\mathcal{L}(X | X>Z^{j})$ if $X_{i}^{j-1}\leq Z^{j}$, and remains unchanged if $X_{i}^{j-1}> Z^{j}$. Notice that the condition $(\sigma^j)^{-1}(i)\leq k$ is equivalent to $i\in\left\{\sigma^{j}(1),\ldots,\sigma^{j}(k)\right\}$.


\item Define $Z^{j+1}=X_{(k)}^{j}$, the $k$-th order statistics of the sample $X^{j}=(X_{1}^{j},\ldots,X_{n}^{j})$, and $\sigma^{j+1}$ the (a.s.) unique associated permutation: $X_{\sigma^{j+1}(1)}^{j}<\ldots<X_{\sigma^{j+1}(n)}^{j}$. 

\item Finally increment $j\leftarrow j+1$.

\end{itemize}

\noindent
{\bf End of the algorithm:}
Define $J^{n,k}(x)=j-1$ as the (random) number of iterations. Notice that $J^{n,k}(x)$ is such that $Z^{J^{n,k}(x)} < a$ and $Z^{J^{n,k}(x)+1} \ge a$.
\end{algo}

The estimator of the probability $P(x)$ is defined by
\begin{equation}\label{eq:estimator}
\hat{p}^{n,k}(x)=C^{n,k}(x)\left(1-\frac{k}{n}\right)^{J^{n,k}(x)},
\end{equation}
with
\begin{equation}\label{eq:corrector}
C^{n,k}(x)=\frac{1}{n}{\rm Card}\left\{i ;\, X_{i}^{J^{n,k}(x)} \ge a\right\}.
\end{equation}

When $x=0$, to simplify notations we set $\hat{p}^{n,k}=\hat{p}^{n,k}(0)$.


\subsection{The Central Limit Theorem}\label{subsect:CLT}

The main result of the paper is the following asymptotic normality result.
\begin{theorem}\label{th}
Under assumption \ref{ass:cdf:c0}, we have the following convergence in law, when $n\rightarrow +\infty$, for any fixed $k$ and $a$:
\begin{equation}
\sqrt{n}\bigl(\hat{p}^{n,k}-p\bigr)\rightarrow \mathcal{N}\bigl(0,-p^2\log(p)\bigr).
\end{equation}
\end{theorem}

Notice that the asymptotic variance does not depend on $k$. This result can be used to define asymptotic normal confidence intervals, for one realization of the algorithm and $n\rightarrow +\infty$. However, the speed of convergence is not known and may depend on the estimated probability $n$ and on the parameter $k$.

Thanks to Theorem \ref{th}, we can study the cost of the use of one realization of the AMS algorithm to obtain a given accuracy when $n\rightarrow +\infty$. In \cite{brehier-lelievre-rousset-2014}, the cost was analyzed when using a sample of $M$ independent realizations of the algorithm, giving an empirical estimator, and the analysis was based on an asymptotic analysis of the variance in the large $n$ limit.

Let $\epsilon$ be some fixed tolerance error, and $\alpha>0$. Denote $r_{\alpha}$ such that $\PP(Z\in [-r_\alpha,r_{\alpha}])=1-\alpha$, where $Z$ is a standard Gaussian random variable.

Then for $n$ large, an asymptotic confidence interval with level $1-\alpha$, centered around $p$, is
\[
\left[p-\frac{r_\alpha\sqrt{-p^2\log(p)}}{\sqrt{n}},p+\frac{r_\alpha\sqrt{-p^2\log(p)}}{\sqrt{n}}\right].
\]
Then the $\epsilon$ error criterion $|\hat{p}^{n,k}-p|\leq \epsilon$ is achieved for $n$ of size $\frac{-p^2\log(p)r_{\alpha}^{2}}{\epsilon^2}$.

However, in average one realization of the AMS algorithm requires a number of steps of the order $-n\log(p)/k$, with $k$ random variables sampled at each iteration. Another cost comes from the sorting of the replicas at initialization, and the insertion at each iteration of the $k$ new sampled replicas in the sorted ensemble of the non-resampled ones. Thus the cost to achieve an accuracy of size $\epsilon$ is in the large $n$ regime of size $n\log(n)\bigl(-p^2\log(p)\bigr)$, not depending on $k$.

This cost must be compared with the one when using a pure Monte-Carlo approximation, with an ensemble of non-interacting replicas of size $n$: thanks to the Central Limit Theorem, the tolerance criterion error $\epsilon$ is satisfied for $n$ of size $\frac{-p(1-p)r_{\alpha}^{2}}{\epsilon^2}$. Despite the $\log(n)$ factor in the AMS case, the performance is improved since $p^2\log(p)=\text{o}(p)$ when $p\rightarrow 0$.

\begin{remark}
Unlike in \cite{brehier-lelievre-rousset-2014}, here we are not able to compare the different choices for $k$, since at the limit $n\rightarrow +\infty$ all the involved quantities do not depend on $k$.
\end{remark}

\section{Proof of the Central Limit Theorem}\label{sect:proof}

The proof is divided into the following steps. First, thanks to Assumption \ref{ass:cdf:c0}, we reduce our study to the case when $X$ is distributed according to the exponential law with parameter $1$: $\PP(X>z)=\exp(-z)$ for any $z>0$. Second, we explain why it is simpler to look at $\log(\hat{p}^{n,k}(x))$ instead of $\hat{p}^{n,k}(x)$, and to use the delta-method to get the result. The third step is the introduction of the characteristic function of $\log(\hat{p}^{n,k}(x))$; following the definition of the algorithm, we prove that it is solution of a functional equation with respect to $x$, which can be transformed into a linear ODE of order $k$. Finally, we analyze the solution of this ODE in the limit $n\rightarrow +\infty$.

\subsection{Reduction to the exponential case}

We start with a reduction of the study to the case when the random variable $X$ is exponentially distributed with parameter $1$. It is based on a change of variable with the following function:
\begin{equation}
\Lambda(x)=-\log\bigl(1-F(x)\bigr).
\end{equation}
It is well-known that $F(X)$ is distributed according to the uniform law on $(0,1)$ (thanks to the continuity Assumption \ref{ass:cdf:c0}), and thus $\Lambda(X)$ is exponentially distributed with parameter $1$. In Section $3$ (see Corollary $3.4$) of \cite{brehier-lelievre-rousset-2014}, it is proved that, the law of the estimator $\hat{p}^{n,k}$ is equal to the one of $\hat{q}^{n,k}$, which is the estimator defined with the same values of the parameters $n$ and $k$, but with two differences: the random variable is exponentially distributed with parameter $1$, and the stopping level is $\Lambda(a)$. Notice that $\mathbb{E}\bigl[\hat{q}^{n,k}\bigr]=\exp\bigl(-\Lambda(a)\bigr)=1-F(a)=p$ (by the unbiasedness result of \cite{brehier-lelievre-rousset-2014}). 

Since the arguments are intricate, we do not repeat them here and we refer the interested reader to \cite{brehier-lelievre-rousset-2014}; from now on, we thus assume:
\begin{assumption}\label{ass:expo}
Assume that $X$ is exponentially distributed with parameter $1$: we denote $\mathcal{L}(X)=\mathcal{E}(1)$.
\end{assumption}

As a consequence, it is enough to prove the following result:
\begin{proposition}\label{propo:CLT}
When $n\rightarrow +\infty$
\begin{equation}
\sqrt{n}\bigl(\hat{p}^{n,k}-p\bigr)\rightarrow \mathcal{N}\bigl(0,a\exp(-2a)\bigr).
\end{equation}
\end{proposition}

We will need the following notations:
\begin{itemize}
\item $f(z)=\exp(-z)\ind_{z>0}$ (resp. $F(z)=\bigl(1-\exp(-z)\bigr)\ind_{z>0}$) for the density (resp. the cumulative distribution function) of the exponential law $\mathcal{E}(1)$ with parameter $1$.
\item $f_{n,k}(z)=k\binom{n}{k}F(z)^{k-1}f(z)\bigl(1-F(z)\bigr)^{n-k}$ for the density of the $k$-th order statistics of a sample of size $n$, made of independent and exponentially distributed random variables, with parameter $1$.
\end{itemize}

Finally, in order to deal with the conditional distributions $\mathcal{L}(X|X>x)$ (under Assumption \ref{ass:expo} it is just a shifted exponential distribution $x+\mathcal{E}(1)$ thanks to the loss of memory property of the exponential law) in the algorithm, we set for any $x\geq 0$ and any $y\geq 0$
\begin{equation}
\begin{gathered}
f(y;x)=f(y-x), \quad F(y;x)=F(y-x),\\
f_{n,k}(y;x)=f_{n,k}(y-x),\\
F_{n,k}(y)=\int_{-\infty}^{y}f_{n,k}(z)dz, \quad F_{n,k}(y;x)=F_{n,k}(y-x).
\end{gathered}
\end{equation}

The following result will be very useful:
\begin{equation}\label{formula:d/dxk}
\left\{
\begin{array}{l}
\begin{gathered}
\frac{d}{dx}f_{n,1}(y;x)=nf_{n,1}(y;x).\\
\text{ for $k \in \{2, \ldots , n-1\}$}, \,\frac{d}{dx}f_{n,k}(y;x)=(n-k+1)\left(f_{n,k}(y;x)-f_{n,k-1}(y;x)\right).
\end{gathered}
\end{array}
\right.
\end{equation} 
The proof follows from straightforward computations.

\subsection{Proof of the CLT} 

The first important idea is to prove the Proposition \ref{propo:CLT} for any initial condition $x\in[0,a]$, not only for $x=0$:
\begin{equation}
\sqrt{n}\bigl(\hat{p}^{n,k}(x)-p(x)\bigr)\rightarrow \mathcal{N}\bigl(0,(a-x)\exp(-2(a-x))\bigr).
\end{equation}

A natural idea is to introduce the characteristic function of $\hat{p}^{n,k}(x)$, and to follow the strategy developed in the previous work. Nevertheless, we are not able to derive a useful functional equation with respect to the $x$ variable. A better strategy is to look at the logarithm of the estimator, and to use a particular case of the delta-method (see for instance \cite{vdv}, Section $3$): if for some sequence of real random variables $(\theta_n)_{n\in\mathbb{N}}$ and some $\theta\in\mathbb{R}$ we have $\sqrt{n}\bigl(\theta_n-\theta)\rightarrow \mathcal{N}(0,\sigma^2)$, then $\sqrt{n}\bigl(\exp(\theta_n)-\exp(\theta)\bigr)\rightarrow \mathcal{N}\bigl(0,\exp(2\theta)\sigma^2\bigr)$ - both in distribution when $n\rightarrow +\infty$.

We thus introduce for any $t\in\R$ and any $0\leq x\leq a$
\begin{equation}\label{eq:phi}
\phi_{n,k}(t,x):=\EE\Bigl[\exp\Bigl(it\sqrt{n}\bigl(\log(\hat{p}^{n,k}(x))-\log(p(x))\bigr)\Bigr)\Bigr].
\end{equation}

For technical reasons, we also set
\begin{equation}\label{eq:chi}
\chi_{n,k}(t,x):=\EE\Bigl[\exp\Bigl(it\sqrt{n}\hat{p}^{n,k}(x)\Bigr)\Bigr]=\exp\bigl(it\sqrt{n}(x-a)\bigr)\phi_{n,k}(t,x),
\end{equation}
on which we derive a functional equation, in order to prove the following result, which implies Proposition \ref{propo:CLT} with the choice $x=0$.

\begin{proposition}\label{propo:chara}
When $n\rightarrow +\infty$, for any $0\leq x\leq a$ and any $t\in\R$
\begin{equation}
\phi_{n,k}(t,x)\rightarrow \exp\left(\frac{t^2(x-a)}{2}\right).
\end{equation}
\end{proposition}

The rest of this section is devoted to the proof of Proposition \ref{propo:chara}, thanks to the following series of results.

\begin{lemma}[Functional Equation]\label{lem:eq_func}
For any $n\in\N$ and any $k\in\left\{1,\ldots,n-1\right\}$, and for any $t\in\R$, the function $x\mapsto\chi_{n,k}(t,x)$ is solution of the following functional equation (with unknown $\chi$): for any $0\leq x\leq a$
\begin{eqnarray}\label{eq:funct_chi}
\chi(t,x)&=&e^{it\sqrt{n}\log(1-\frac{k}{n})}\int_{x}^{a}\chi(t,y)f_{n,k}(y;x)\, \mathrm dy \\
&&+  \sum_{l=0}^{k-1}e^{it\sqrt{n}\log(1-\frac{l}{n})}\mathbb{P}(S(x)_{(l)}^{n} < a \leq S(x)_{(l+1)}^{n}),
\end{eqnarray}
where $(S(x)_{j}^{n})_{1 \leq j \leq n}$ are iid with law $\mathcal{L}(X|X>x)$ and where $S(x)_{(l)}^{n}$ is the $l$-th order statistics of this sample (with convention $S(x)_{(0)}^{n}=x$).
\end{lemma}

\begin{remark}
The function $\mapsto\phi_{n,k}(t,x)$ is solution of a more complex functional equation: for any $0\leq x\leq a$
\begin{eqnarray*}
\phi_{n,k}(t;x)&=&\int_{x}^{a}e^{it\sqrt{n}\log(1-k/n)}e^{-it\sqrt{n}(x-y)}\phi_{n,k}(t;y)f_{n,k}(y;x)\mathrm dy\\
&&+\sum_{l=0}^{k-1}e^{-it\sqrt{n}(x-a)}e^{it\sqrt{n}\log(1-\frac{l}{n})}\mathbb{P}(S(x)_{(l)}^{n} < a \leq S(x)_{(l+1)}^{n}).
\end{eqnarray*}
The advantage of the equation on $\chi_{n,k}$ is that we are able to deduce that $x\mapsto\chi_{n,k}(t,x)$ is solution of a linear ODE with constant coefficients. On $x\mapsto\phi_{n,k}(t,x)$ the ODE would not have constant coefficients and would be more difficult to express.
\end{remark}

\begin{remark}
If instead of considering $\log(\hat{p}^{n,k})(x)$ we try to prove directly the central limit theorem on $\hat{p}^{n,k}(x)$, the same strategy leads to a functional equation where on the left-hand side we have the parameter $t$, and on the right-hand side, in the integrand we have the parameter $t\bigl(1-k/n\bigr)$.
\end{remark}

\Proof
The idea (like in the proof of Proposition $4.2$ in \cite{brehier-lelievre-rousset-2014}) is to decompose the expectation according to the value of the first level $Z^1=X_{(k)}^{0}$. On the event $\left\{Z^1>a\right\}=\left\{J^{n,k}(x)=0\right\}$, the algorithm stops and $\hat{p}^{n,k}(x)=\frac{n-l}{n}$ for the unique $l\in\left\{0,\ldots,k-1\right\}$ such that $S(x)_{(l)}^{n} < a \leq S(x)_{(l+1)}^{n}$. Thus
\begin{equation}
\mathbb{E}[e^{it\sqrt{n}\log(\hat{p}^{n,k}(x))} \ind_{J^{n,k}(x)=0}] = \sum_{l=0}^{k-1}e^{it\sqrt{n}\log(1-\frac{l}{n})} \mathbb{P}(S(x)_{(l)}^{n} < a \leq S(x)_{(l+1)}^{n}).
\end{equation}
If $Z^1<a$, for the next iteration the algorithm restarts at the initial condition $Z^1$, and
\begin{equation}
\begin{aligned}
\footnotesize{\mathbb{E}[e^{it\sqrt{n}\log(\hat{p}^{n,k}(x))}}&\footnotesize{ \ind_{J^{n,k}(x)>0}]}\\
&=\footnotesize{\mathbb{E}\left[e^{it\sqrt{n}\log(1 - \frac{k}{n})}\mathbb{E}[e^{it\sqrt{n}\log\left(C^{n,k}(x)(1 - \frac{k}{n})^{J^{n,k}(x)-1}\right)}| Z^1]\ind_{Z^1<a} \right]} \\
&=\footnotesize{e^{it\sqrt{n}\log(1 - \frac{k}{n})}\mathbb{E}\left[\mathbb{E}[e^{it\sqrt{n}\log\left(\hat{p}^{n,k}(Z^1)\right)}| Z^1]\ind_{Z^1<a} \right]} \\
&=\footnotesize{e^{it\sqrt{n}\log(1 - \frac{k}{n})}\mathbb{E}\left[\chi_{n,k}(t,Z^1)\ind_{Z^1<a} \right]} \\
&=\footnotesize{e^{it\sqrt{n}\log(1 - \frac{k}{n})}\int_{x}^{a}\chi_{n,k}(t,y)f_{n,k}(y;x)\, \mathrm dy} .
\end{aligned}
\end{equation}
With the two identities and the definition \eqref{eq:chi} of $\chi_{n,k}$, we thus obtain \eqref{eq:funct_chi}.
\EProof

We exploit the functional equation \eqref{eq:funct_chi} on $x\mapsto\chi_{n,k}(t,x)$, to prove that this function is solution of a linear ODE with constant coefficients.
\begin{lemma}[ODE]\label{lem:ODE}
Let $n$ and $k \in \{1, \ldots , n-2\}$ be fixed. There exist real numbers $\mu^{n,k}$ and $(r_{m}^{n,k})_{0\leq m\leq k-1}$, depending only on $n$ and $k$, such that for any fixed $t\in \R$, the function $x\mapsto\chi_{n,k}(t,x)$ satisfy the following Linear Ordinary Differential Equation (ODE) of order $k$: for $x \in [0,a]$
\begin{equation}\label{eq:ODE}
\frac{d^k}{dx^k}\chi_{n,k}(t,x)=e^{it\sqrt{n}\log(1-\frac{k}{n})}\mu^{n,k}\chi_{n,k}(t,x)+\sum_{m=0}^{k-1}r_{m}^{n,k}\frac{d^m}{dx^m}\chi_{n,k}(t,x).
\end{equation}
The coefficients $\mu^{n,k}$ and $(r_{m}^{n,k})_{0\leq m\leq k-1}$ satisfy the following properties:
\begin{equation}\label{eq:rec_coeffs}
\begin{gathered}
\mu^{n,k}=(-1)^{k}n\ldots (n-k+1)\\
\lambda^k-\sum_{m=0}^{k-1}r_{m}^{n,k} \, \lambda^m=(\lambda-n)\ldots (\lambda-n+k-1) \quad \text{ for all } \lambda\in\R.
\end{gathered}
\end{equation}
\end{lemma}

\Proof
The proof follows the same lines as Proposition $6.4$ in \cite{brehier-lelievre-rousset-2014}. We introduce the following notation
$$\Theta_{n,k}(t,x):=\sum_{l=0}^{k-1}e^{it\sqrt{n}\log(1-\frac{l}{n})}\mathbb{P}(S(x)_{(l)}^{n} < a \leq S(x)_{(l+1)}^{n}).$$
Then by recursion, we can prove that for $0\leq l\leq k-1$ and for any $x\leq a$ and $t\in\R$
\begin{eqnarray}\label{eq:recur_p}
\frac{d^l}{dx^l}\left(\chi_{n,k}(t,x)-\Theta_{n,k}(t,x)\right)&=&\mu_{l}^{n,k}e^{it\sqrt{n}\log(1-\frac{k}{n})}\int_{x}^{a}\chi_{n,k}(t,y)f_{n,k-l}(y;x)\, \mathrm dy\nonumber\\
&&+\sum_{m=0}^{l-1}r_{m,l}^{n,k}\frac{d^m}{dx^m}\left(\chi_{n,k}(t,x)-\Theta_{n,k}(t,x)\right).
\end{eqnarray}
The idea for the recursion is to differentiate \eqref{eq:recur_p} with respect to $x$ and to use \eqref{formula:d/dxk}. The coefficients satisfy the following recursion (notice that they do not depend on $t$):
\begin{equation}\label{eq:recursion}
\begin{gathered}
\mu_{0}^{n,k}=1,\mu_{l+1}^{n,k}=-(n-k+l+1)\mu_{l}^{n,k};\\
\begin{cases}
r_{0,l+1}^{n,k}=-(n-k+l+1)r_{0,l}^{n,k}, \quad \text{if } l>0,\\
r_{m,l+1}^{n,k}=r_{m-1,l}^{n,k}-(n-k+l+1)r_{m,l}^{n,k}, \quad 1\leq m\leq l,\\
r_{l,l}^{n,k}=-1.
\end{cases}
\end{gathered}
\end{equation}
Differentiating once more when $l=k-1$, we obtain
\begin{eqnarray}\label{eq:recur_p_fin}
\frac{d^k}{dx^k}\left(\chi_{n,k}(t,x)-\Theta_{n,k}(t,x)\right)&=&\mu^{n,k}e^{it\sqrt{n}\log(1-\frac{k}{n})}\chi_{n,k}(t,x)\nonumber\\
&&+\sum_{m=0}^{k-1}r_{m}^{n,k}\frac{d^m}{dx^m}\left(\chi_{n,k}(t,x)-\Theta_{n,k}(t,x)\right),
\end{eqnarray}
with $\mu^{n,k}:=\mu_{k}^{n,k}$ and $r_{m}^{n,k}:=r_{m,k}^{n,k}$. These coefficients are the same as in \cite{brehier-lelievre-rousset-2014}; in particular it was proved there that the claimed properties are satisfied (using the unbiasedness result for the estimation of the probability, see Section $6.4$ there). The polynomial equality of \eqref{eq:rec_coeffs} is in fact equivalent to the following identity: for all $j\in\left\{0,\ldots,k-1\right\}$
$$\frac{d^k}{dx^k}\exp\left((n-k+j+1)(x-a)\right)=\sum_{m=0}^{k-1}r_{m}^{n,k}\frac{d^m}{dx^m}\exp\left((n-k+j+1)(x-a)\right).$$
Using the expression of the order statistics in the exponential case,  $\Theta_{n,k}(t,.)$ is a linear combination of the exponential functions $z\mapsto \exp(-nz),\ldots,\exp(-(n-k+1)z)$; therefore
$$\frac{d^k}{dx^k}\Theta_{n,k}(t,x)=\sum_{m=0}^{k-1}r_{m}^{n,k}\frac{d^m}{dx^m}\Theta_{n,k}(t,x),$$
and thus \eqref{eq:recur_p_fin} gives \eqref{eq:ODE}.
\EProof

To conclude, it remains to express the solution of \eqref{eq:ODE}, and to analyze the asymptotics $n\rightarrow +\infty$. Since the ODE is of order $k$, in order to uniquely determine the solution, more information is required: for instance we need to know the vector of the derivatives of order $0,1,\ldots , k-1$ of $x\mapsto\chi_{n,k}(t,x)$ at some point. From the functional equation \eqref{eq:funct_chi}, we are able to show the following asymptotic result at point $a$, when $n\rightarrow +\infty$.
\begin{lemma}[Initial condition]\label{lem:init_cond}
For any fixed $k\in\left\{1,\ldots,\right\}$ and any $t\in\R$, we have
\begin{equation}\label{eq:init_cond}
\left\{
\begin{array}{l}
\chi_{n,k}(t,a) = 1 \\
\frac{d^{m}}{dx^{m}}\chi_{n,k}(t,x)\Big|_{x=a} \underset{n \to \infty} = O(\frac{1}{\sqrt{n}})n^{m} \quad \text{ if } m\in\left\{1,\ldots,k-1\right\}.
\end{array}
\right.
\end{equation}
\end{lemma}

\Proof
The equality $\chi_{n,k}(t,a) = 1$ is trivial, since $\hat{p}^{n,k}(a)=1$. From \eqref{eq:recur_p} and \eqref{eq:recur_p_fin}, we immediately get (by recursion) that for $1\leq m\leq k-1$
$$\frac{d^{m}}{dx^{m}}\chi_{n,k}(t,x)\Big|_{x=a}=\frac{d^{m}}{dx^{m}}\Theta_{n,k}(t,x)\Big|_{x=a}.$$
We introduce the following decomposition
\begin{align*}
\Theta_{n,k}(t,x)&=\sum_{l=0}^{k-1}(e^{it\sqrt{n}\log(1-\frac{l}{n})}\mathbb{P}(S(x)_{(l)}^{n} < a \leq S(x)_{(l+1)}^{n})\\
&=\sum_{l=0}^{k-1}\left(e^{it\sqrt{n}\log(1-\frac{l}{n})}-1\right)\bigl(F_{n,l}(a;x)-F_{n,l+1}(a;x) \bigr)\\
&~+\sum_{l=0}^{k-1}\mathbb{P}(S(x)_{(l)}^{n} < a \leq S(x)_{(l+1)}^{n})\\
&=:\Omega_{n,k}(t,x)+1-F_{n,k}(a;x),
\end{align*}
where $F_{n,l}$ denotes the cumulative distribution function of the $l$-th order statistics (with the convention $F_{n,0}(a;x)=1$ for $x\leq a$).

Thanks to \eqref{formula:d/dxk} and a simple recursion on $l$, it easy to prove that for any $0\leq l\leq k$ and any $m\geq 1$
\begin{equation}\label{eq:rec_init}
\frac{d^{m}}{dx^{m}}F_{n,l}(a;x)\Big|_{x=a}=\text{O}(n^m),
\end{equation}
which immediately yields 
$$\frac{d^{m}}{dx^{m}}\Omega_{n,k}(t,x)\Big|_{x=a} \underset{n \to \infty} = O(\frac{1}{\sqrt{n}})n^{m}.$$
In fact, it is possible to prove a stronger result: if $1\leq l\leq k$ and $0\leq m<l$ then
$$\frac{d^{m}}{dx^{m}}F_{n,l}(a;x)\Big|_{x=a}=0,$$
by recursion on $l$ and using \eqref{formula:d/dxk} recursively on $m$. We thus obtain for $1\leq m\leq k-1$
$$\frac{d^{m}}{dx^{m}}\bigl(1-F_{n,k}(a;x)\bigr)\Big|_{x=a}=0.$$
This concludes the proof of the Lemma.
\EProof

The last useful result is the following:
\begin{lemma}[Asymptotic expansion]\label{lem:express}
Let $k\in\left\{1,\ldots,\right\}$ and $t\in\R$ be fixed. Then for $n$ large enough, we have
\begin{equation}\label{eq:sol_ODE}
\chi_{n,k}(t,x) = \sum_{l=1}^{k} \eta_{n,k}^{l}(t) e^{\lambda_{n,k}^{l}(t)\left(x-a\right)},
\end{equation}
for complex coefficients satisfying:
\begin{equation}
\left\{
\begin{array}{l}
  \lambda_{n,k}^{1}(t) \underset{n \to \infty} = it\sqrt{n} + \frac{t^{2}}{2} + o(1) ,\\
  \eta_{n,k}^{1}(t) \underset{n \to \infty} \rightarrow 1;
\end{array}
\right.
\end{equation}
and for $2\leq l\leq k$
\begin{equation}
\left\{
\begin{array}{l}
  \lambda_{n,k}^{l}(t) \underset{n \to \infty} \sim n(1-e^{\frac{i2\pi(l-1)}{k}}) ,\\
  \eta_{n,k}^{l}(t) \underset{n \to \infty} \rightarrow 0.
\end{array}
\right.
\end{equation}
\end{lemma}

\Proof
We denote by $(\lambda_{n,k}^{l}(t))_{1\leq l\leq k}$ the roots of the characteristic equation (with unknown $\lambda\in \C$):
$$\frac{(n-\lambda)...(n-k+1-\lambda)}{n...(n-k+1)}-e^{it\sqrt{n}\log(1-\frac{k}{n})}=0$$
By the continuity property of the roots of a complex polynomial of degree $k$ with respect to its coefficients, we have
$$\overline{\lambda}_{n,k}^{l}(t):=\frac{\lambda_{n,k}^{l}(t)}{n} \underset{n \to \infty} \rightarrow\overline{\lambda}_{\infty}^{l},$$
where $(\overline{\lambda}_{\infty}^{l}(t))_{1\leq l\leq k}$ are the roots of $(1-\overline{\lambda})^k=1$: thus, we get $\lambda_{n,k}^{1}(t) \underset{n \to \infty}=\text{o}(n)$,
$$\lambda_{n,k}^{l}(t) \underset{n \to \infty} \sim n(1-e^{\frac{i2\pi(l-1)}{k}}).$$
To study more precisely the asymptotic behavior of $\lambda_{n,k}^{1}(t)$, we postulate an ansatz
$$ \lambda_{n,k}^{1}(t) \underset{n \to \infty} = c_{t}\sqrt{n}+d_{t}+o(1),$$
and identify the coefficients $c_t=it$ and $d_t=t^2/2$ thanks to
\begin{gather*}
\left(1-\frac{c_{t}}{\sqrt{n}} - \frac{d_{t}}{n} + o(1) \right)^{k} \underset{n \to \infty} = 1 - \frac{c_tk}{\sqrt{n}} - \frac{d_tk-\binom{k}{2}c_t^2}{n} + \text{o}\left(\frac{1}{n}\right)\\
e^{it\sqrt{n}\log(1-\frac{k}{n})}\underset{n \to \infty} =1 - \frac{itk}{\sqrt{n}} - \frac{t^{2}k^{2}}{2n} + \text{o}\left(\frac{1}{n}\right).
\end{gather*}
In particular, for $n$ large enough, $(\lambda_{n,k}^{l}(t))_{1\leq l\leq k}$ are pairwise distinct, and \eqref{eq:sol_ODE} follows.

The coefficients $(\eta_{n,k}^{l}(t))_{1\leq l\leq k}$ are then solutions of the following linear system of order $k$:
\begin{equation}
\left\{
\begin{array}{l}
\eta_{n,k}^{1}(t) + ... + \eta_{n,k}^{k}(t) = \chi_{n,k}(a),\\
\eta_{n,k}^{1}(t)\overline{\lambda}_{n,k}^{1}(t) + ... + \eta_{n,k}^{k}(t)\overline{\lambda}_{n,k}^{k}(t) =\frac{1}{n}\frac{d}{dx}\chi_{n,k}(t,a),\\
... \\
\eta_{n,k}^{1}(t)\left(\overline{\lambda}_{n,k}^{1}(t)\right)^{k-1} + ... + \eta_{n,k}^{k}(t)\left(\overline{\lambda}_{n,k}^{k}(t)\right)^{k-1} =\frac{1}{n^{k-1}}\frac{d^{k-1}}{dx^{k-1}}\chi_{n,k}(t,a).
\end{array}
\right.
\end{equation}

Using Cramer's rule, we express each $\eta_{n,k}^{l}(t)$ as a ratio of determinants (the denominator being a Vandermonde determinant is non zero when $n$ is large enough). For $l\in\left\{2,\ldots,k\right\}$, we have
$$\eta_{n,k}^{l}(t)=\frac{\det(M_{n,k}^{l}(t))}{V\Bigl(\overline{\lambda}_{n,k}^{1}(t),\ldots,\overline{\lambda}_{n,k}^{k}(t)\Bigr)}\underset{n\rightarrow +\infty}\rightarrow 0$$
where
$$M_{n,k}^{l}(t) =  \begin{pmatrix}
  1 &  1 & \ldots & 1 & \ldots & 1\\
   \overline{\lambda}_{n,k}^{1}(t) & \overline{\lambda}_{n,k}^{2}(t) & \ldots & O(\frac{1}{\sqrt{n}}) & \ldots & \overline{\lambda}_{n,k}^{k}(t)\\
   \vdots & \vdots & \vdots & \vdots & \vdots & \vdots\\
   \bigl(\overline{\lambda}_{n,k}^{1}(t)\bigr)^{k-1} & \bigl(\overline{\lambda}_{n,k}^{2}(t)\bigr)^{k-1} & \ldots & O(\frac{1}{\sqrt{n}}) & \ldots & \bigl(\overline{\lambda}_{n,k}^{k}(t)\bigr)^{k-1}
\end{pmatrix}$$
is such that $\det(M_{n,k}^{l}(t))\underset{n\rightarrow +\infty}\rightarrow 0$ (since $\overline{\lambda}_{n,k}^{1}(t)\rightarrow 0$), and the denominator is given by a Vandermonde determinant
$
V\Bigl(\overline{\lambda}_{n,k}^{1}(t),\ldots,\overline{\lambda}_{n,k}^{k}(t)\Bigr)\underset{n\rightarrow +\infty}\rightarrow V\Bigl(\overline{\lambda}_{\infty}^{1}(t),\ldots,\overline{\lambda}_{\infty}^{k}(t)\Bigr)\neq 0.
$

Finally, $\eta_{n,k}^{1}(t)=1-\sum_{l=2}^{k}\eta_{n,k}^{l}(t)\underset{n\rightarrow +\infty}\rightarrow 1$.
\EProof

The proof of Proposition \ref{propo:chara} is now straightforward, recalling from \eqref{eq:phi} and \eqref{eq:chi} the relation $\phi_{n,k}(t,x)=\exp\bigl(-it\sqrt{n}(x-a)\bigr)\chi_{n,k}(t,x)$ between the characteristic function $\phi_{n,k}$ and the auxiliary function $\chi_{n,k}$, and taking the limit $n\rightarrow +\infty$ thanks to Lemma \ref{lem:express}.

\section{Numerical results}\label{sect:num}

In this Section, we provide some numerical illustrations of the Central Limit Theorem \ref{th}. We apply the algorithm in the idealized setting with an exponentially distributed random variable with parameter $1$; indeed, we have seen that in the idealized setting, in law the algorithm behaves like this, if the threshold is suitably chosen.

In the simulations we present, the probability we estimate is $\exp(-6)$, which is approximately $2.10^{-3}$.

\begin{figure}[h]
\centering

   \caption{Histograms for $k=10$ and $p=\exp(-6)$: $n=10^2,10^3,10^4$ from left to right}
\subfigure{
	\includegraphics[scale=0.16]{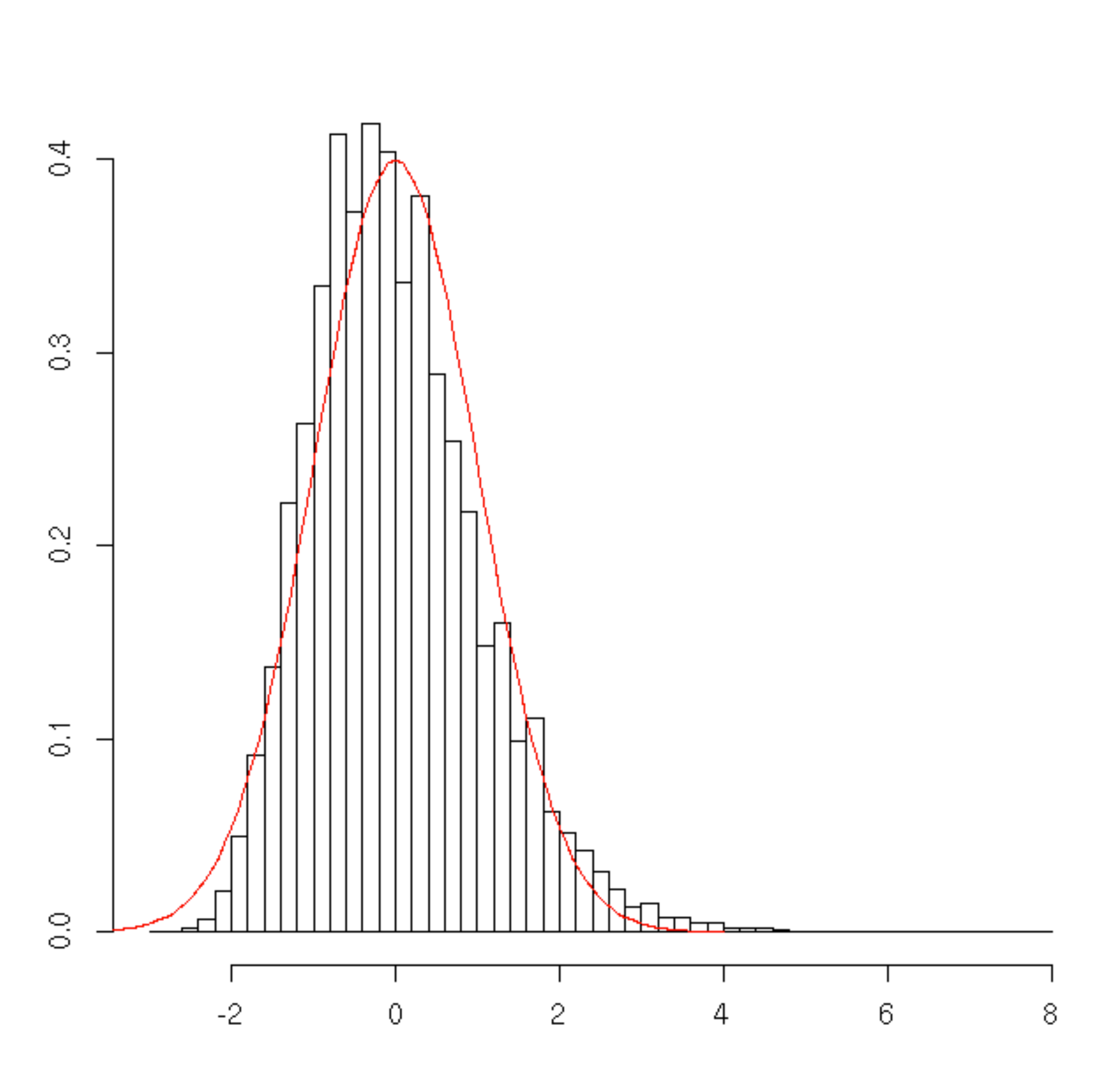}}
\subfigure{
	\includegraphics[scale=0.16]{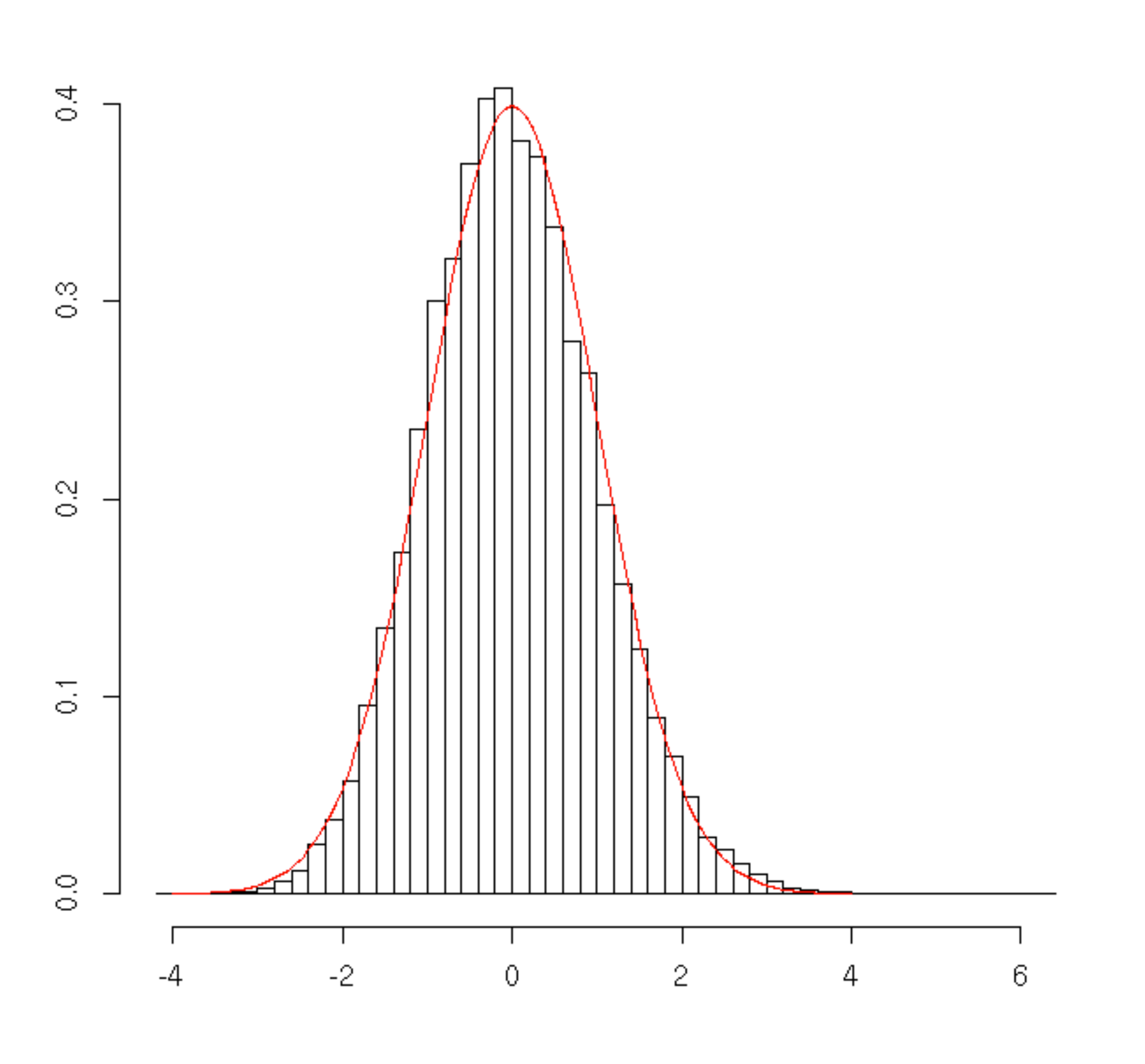}}
\subfigure{
	\includegraphics[scale=0.16]{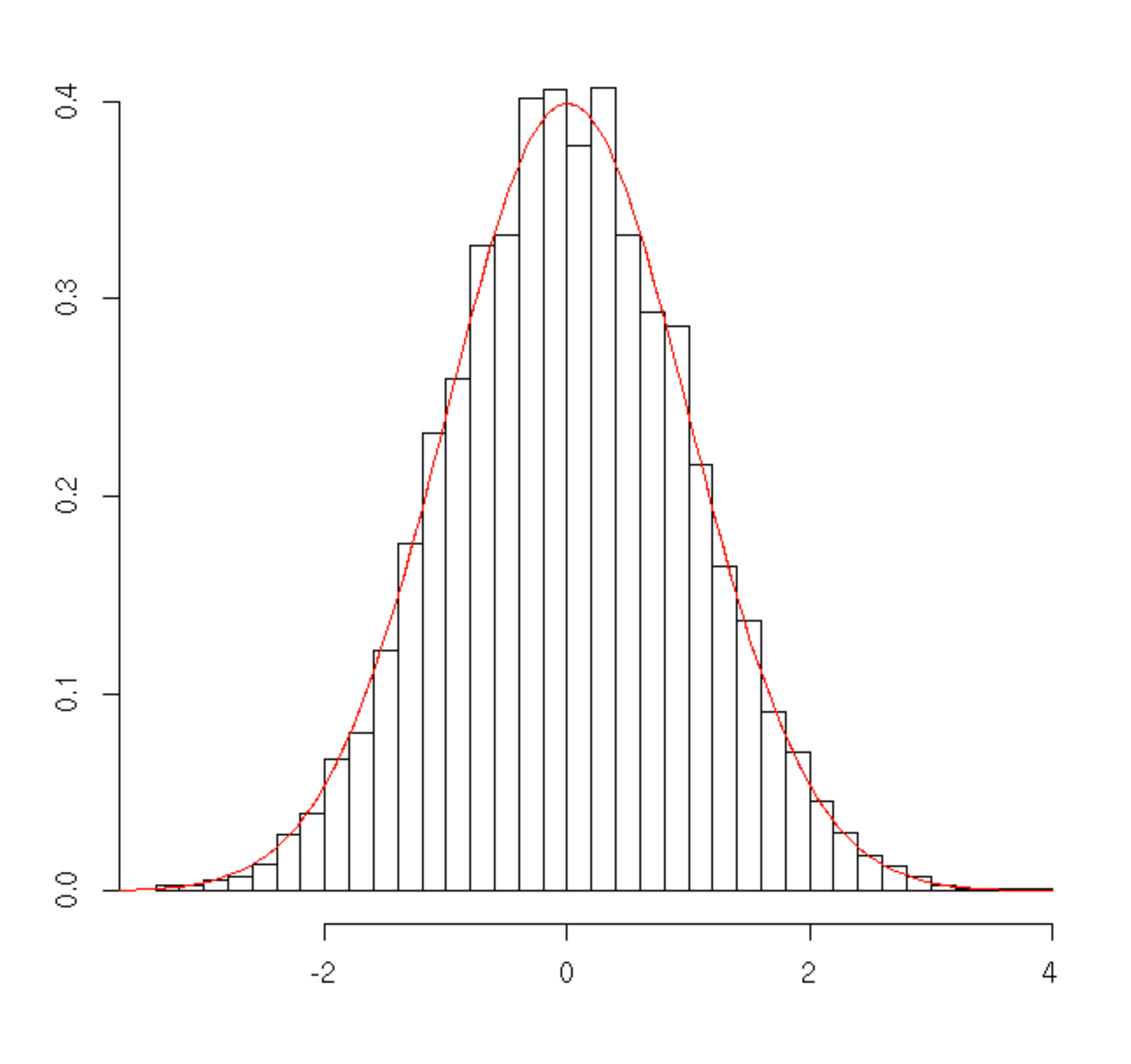}}
\label{fig:n's}
\end{figure}

In figure \ref{fig:n's}, we fix the value $k=10$, and we show histograms for $n=10^2,10^3,10^4$, with different values for the number $M$ independent realizations of the algorithm, such that $nM=10^8$ (we thus have empirical variance of the same order for all cases). When $n$ increases, the normality of the estimator is confirmed. In figure \ref{fig:n's}, we give Q-Q plots, where the empirical quantiles of the sample are compared with the exact quantiles of the standard Gaussian random variable (after normalization).

\begin{figure}[h]
\centering
   \caption{Q-Q plot for $k=10$ and $p=\exp(-6)$: $n=10^2,10^3,10^4$ from left to right}
\subfigure{
	\includegraphics[scale=0.14]{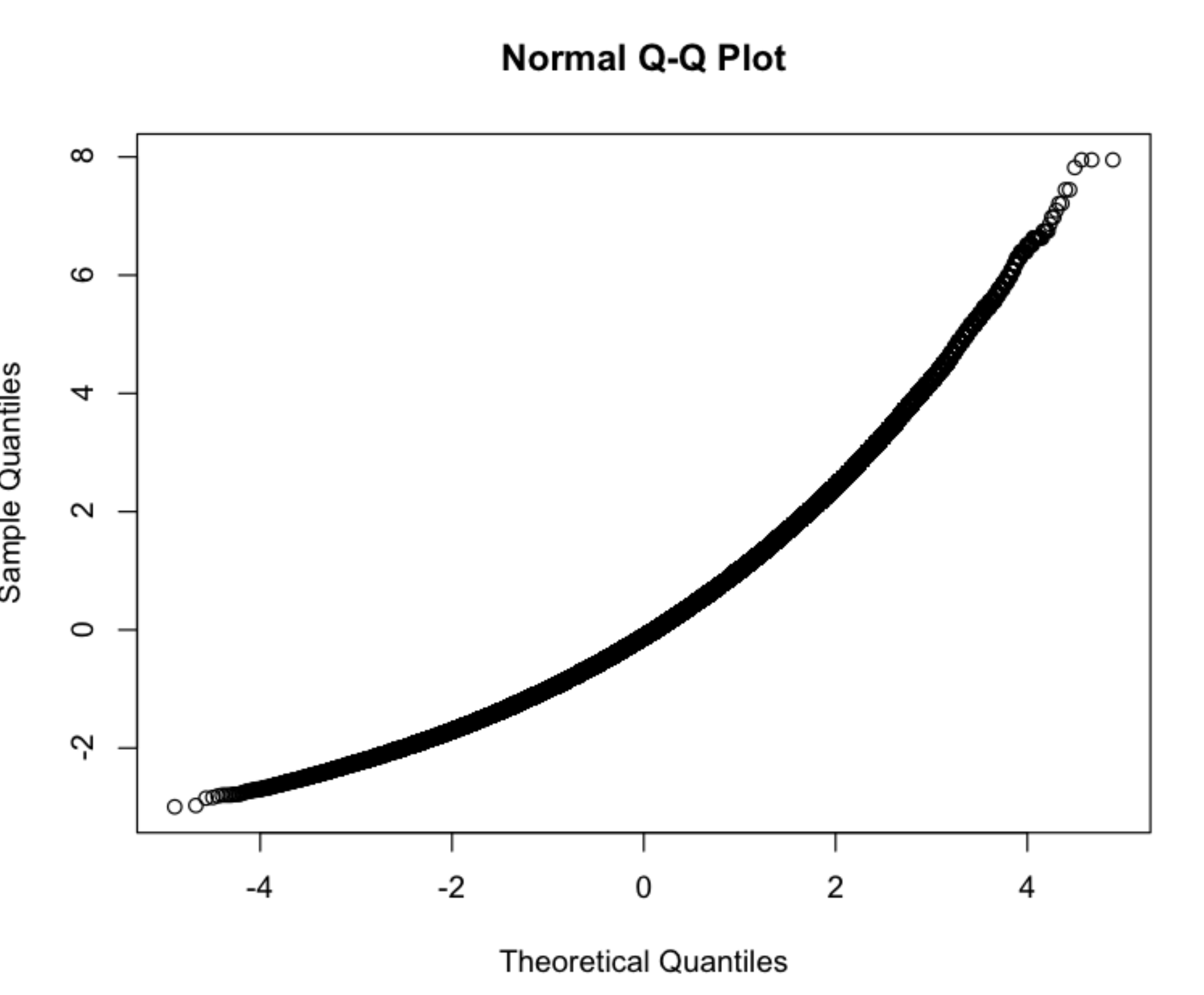}}
\subfigure{
	\includegraphics[scale=0.14]{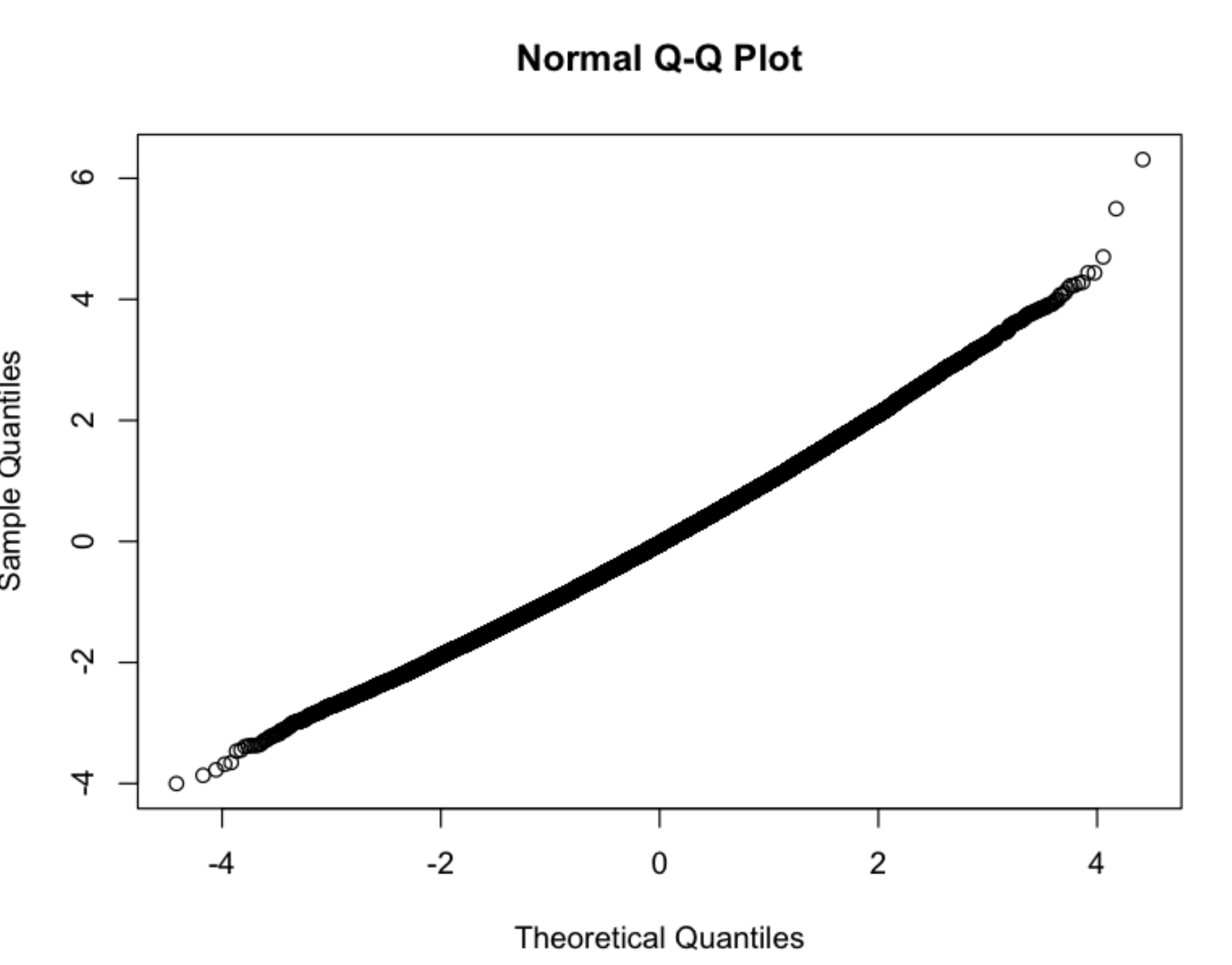}}
\subfigure{
	\includegraphics[scale=0.14]{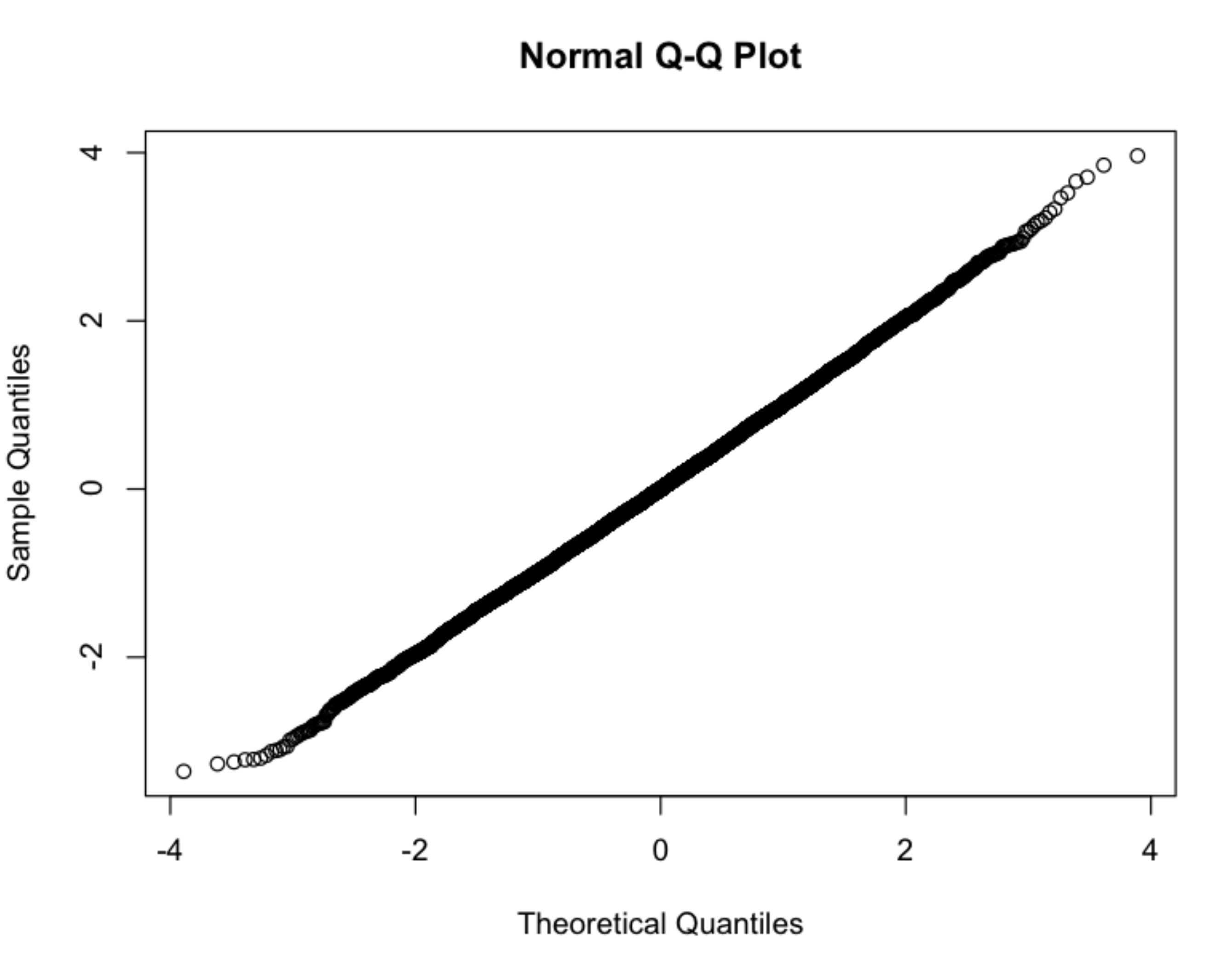}}
\label{fig:qq-n's}
\end{figure}

\begin{figure}[h]
\centering

   \caption{Histograms for $n=10^4$ and $p=\exp(-6)$: $k=1,10,100$ from left to right}
\subfigure{
	\includegraphics[scale=0.16]{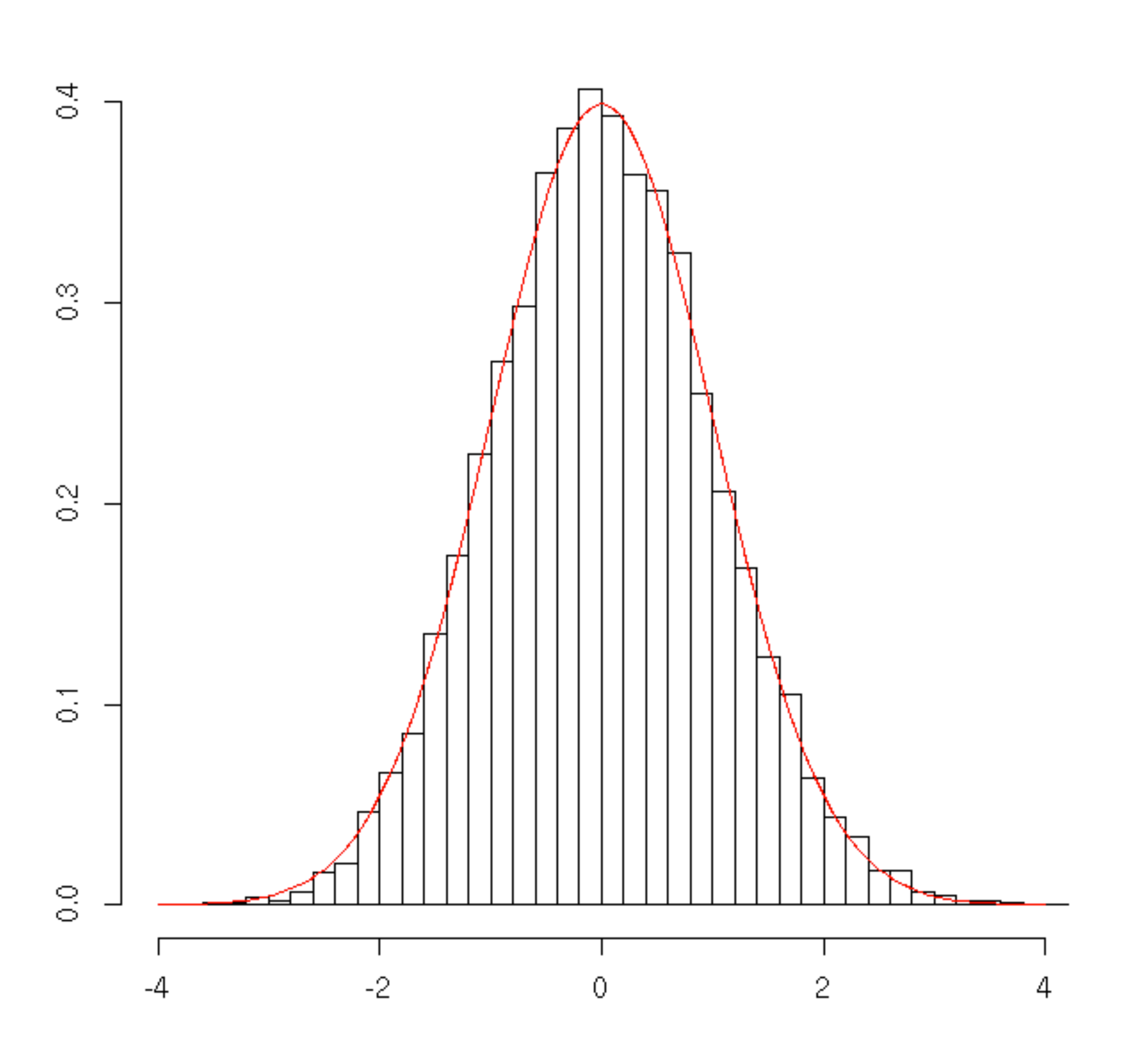}}
\subfigure{
	\includegraphics[scale=0.16]{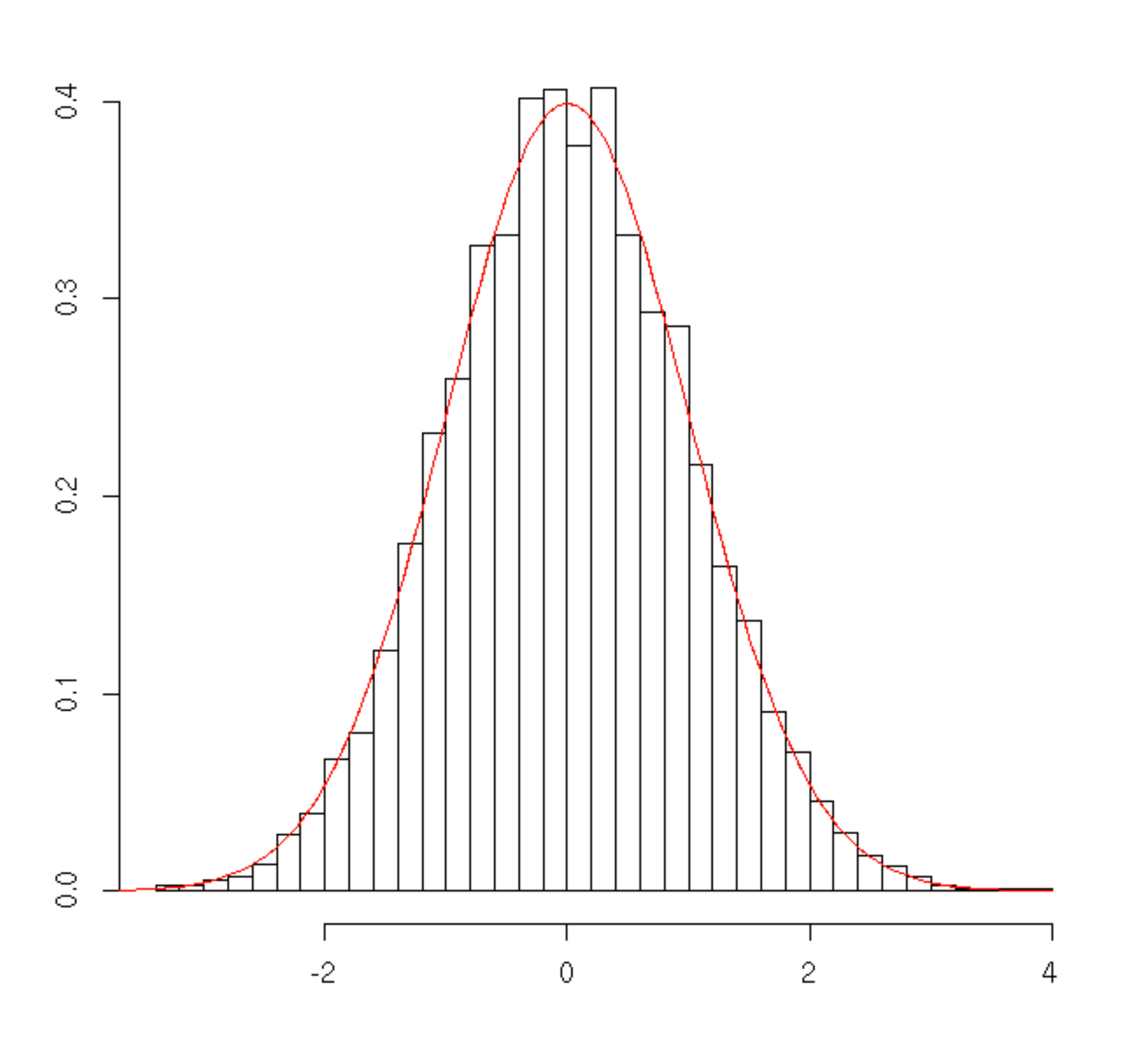}}
\subfigure{
	\includegraphics[scale=0.16]{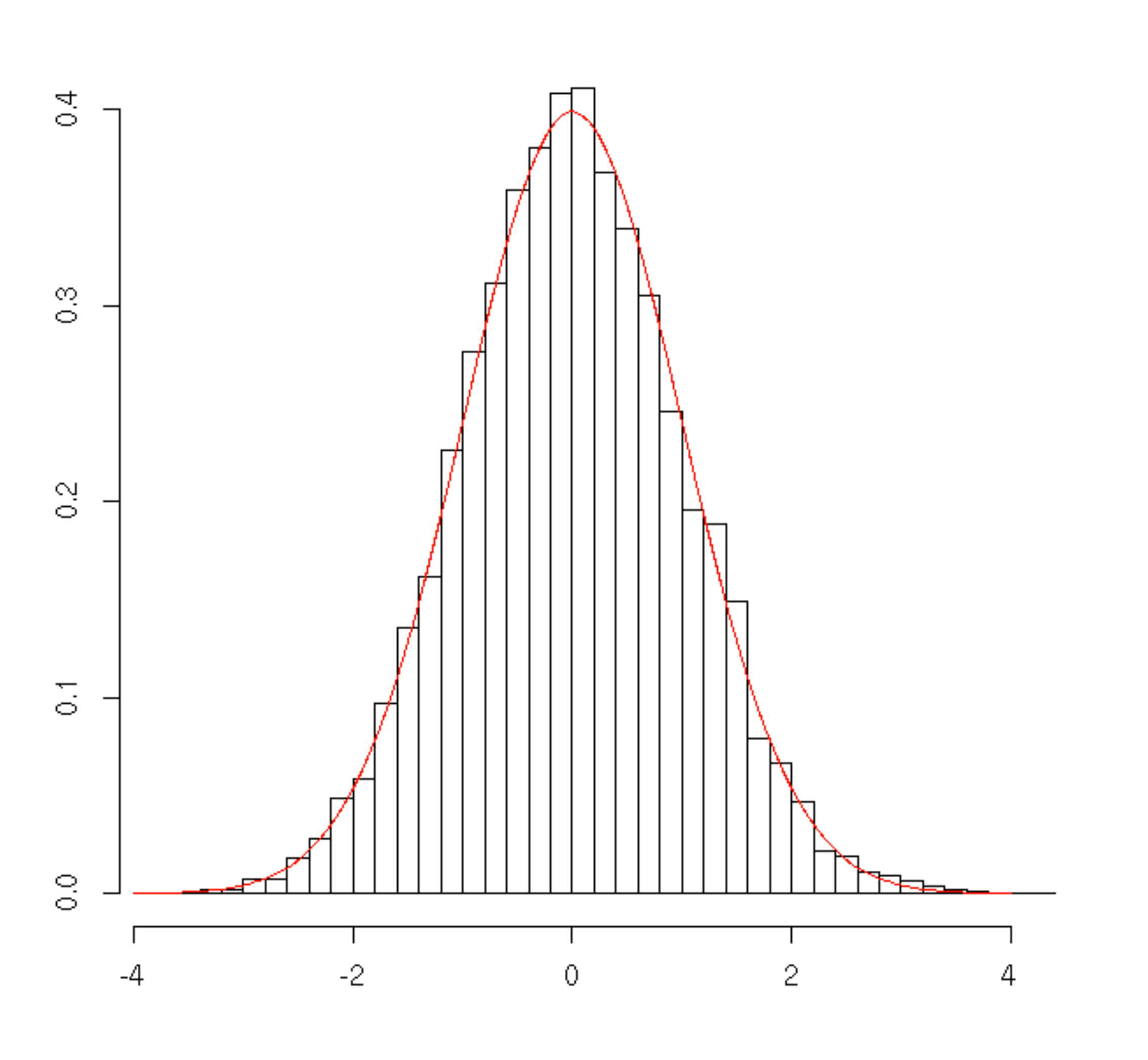}}
\label{fig:k's}
\end{figure}

\begin{figure}[h!]
\centering
   \caption{Q-Q plot for $n=10^4$ and $p=\exp(-6)$: $k=1,10,100$ from left to right}
\subfigure{
	\includegraphics[scale=0.14]{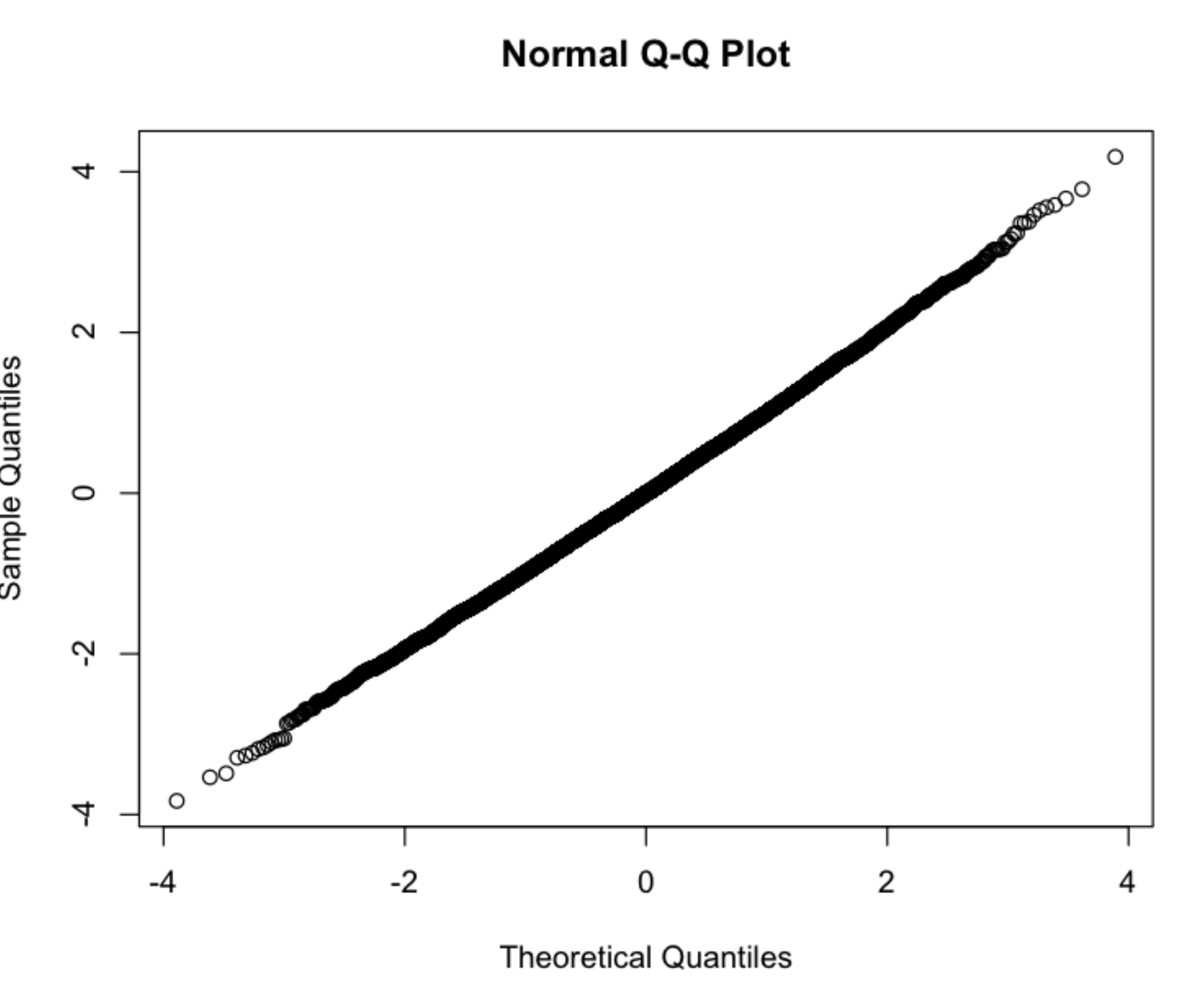}}
\subfigure{
	\includegraphics[scale=0.14]{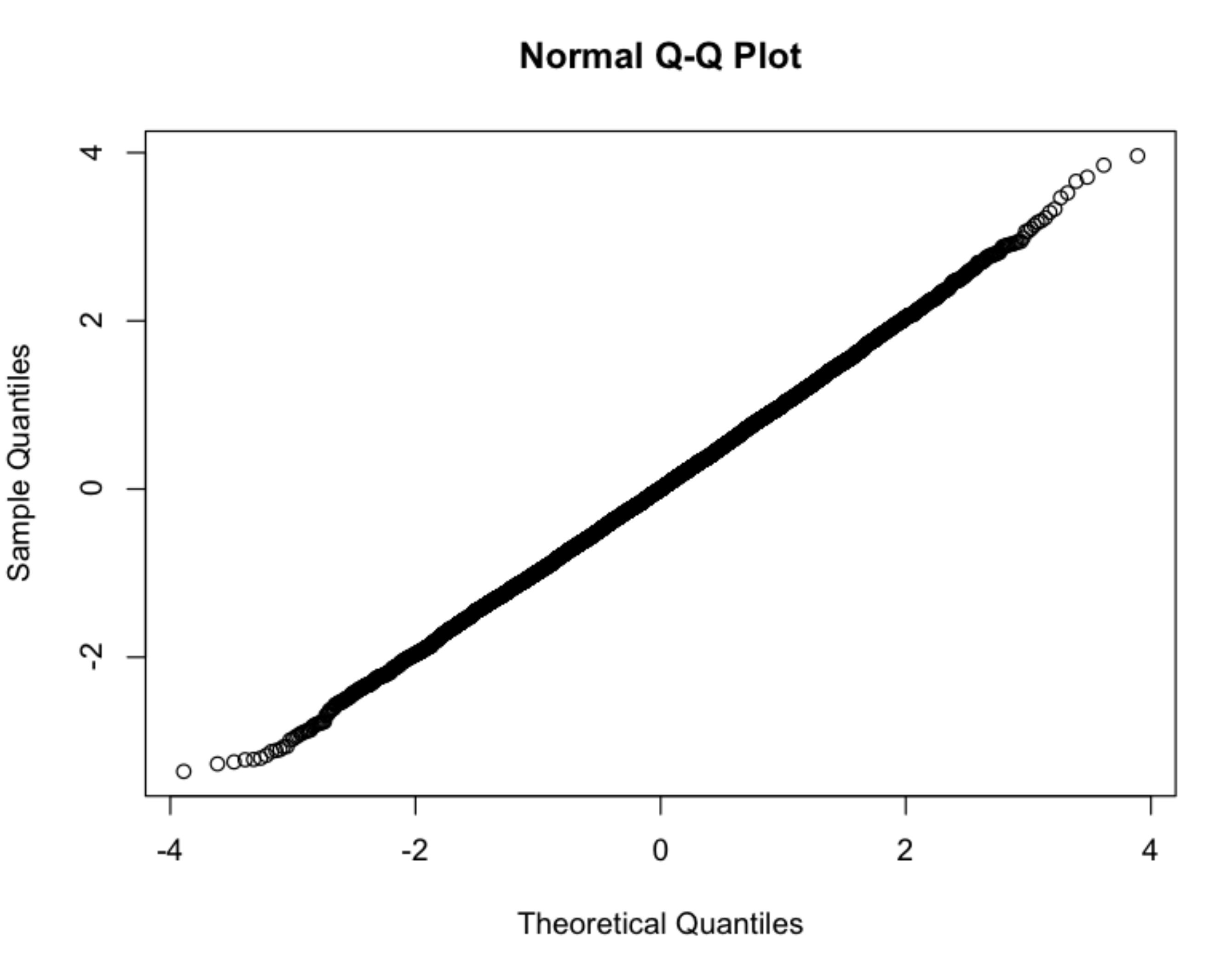}}
\subfigure{
	\includegraphics[scale=0.14]{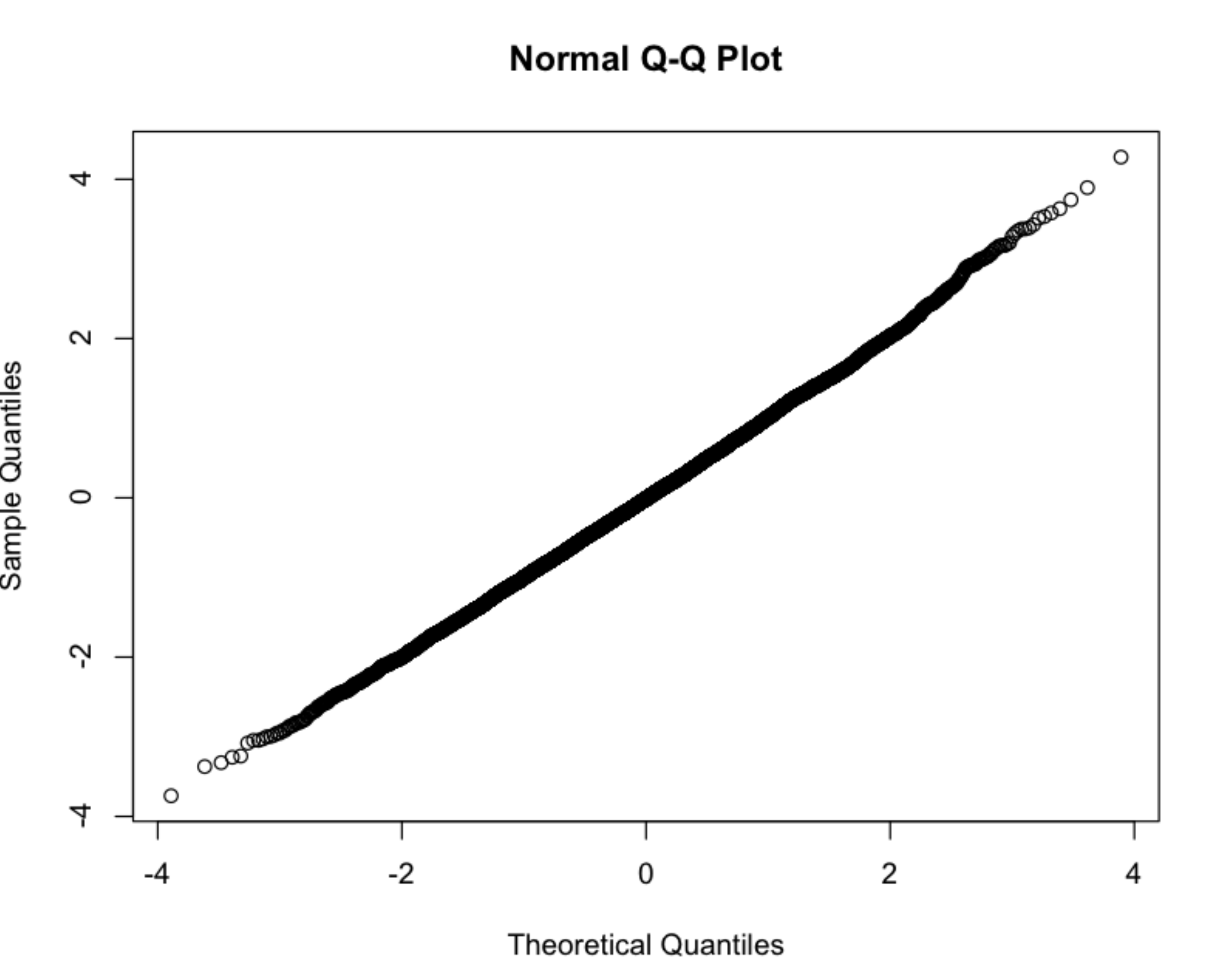}}
\label{fig:qq-k's}
\end{figure}

In Figure \ref{fig:k's}, we show histograms for $M=10^4$ independent realizations of the AMS algorithm with $n=10^4$ and $k\in\left\{1,10,100\right\}$; we also provide associated Q-Q plots in figure \ref{fig:qq-k's}, which prove that in this regime the normality assumption seems to be satisfied for all values of $k$.

\acknowledgements{
C.-E.~B. would like to thank G.~Samaey, T.~Leli\`evre and M.~Rousset for the invitation to give a talk on the topic of this paper at the $11$-th MCQMC Conference, in the special session on Mathematical aspects of Monte Carlo methods for molecular dynamics.
}

\bibliographystyle{spmpsci}
%

\end{document}